\documentstyle[12pt]{article}

\newcommand{\BOX}{\hbox {$\sqcap$ \kern -1em $\sqcup$}}
\newcommand{\qed}{\hskip 3em \hbox{\BOX} \vskip 2ex}
         \newcommand{\Sl}{{\bf sl}}
          \newcommand{\g}{{\bf g}}

        \newtheorem{theorem}{Theorem}
        \newcommand{\et}{\hspace{-0.08in}{\bf .}\hspace{0.1in}}
        \newcommand{\maps}{\colon}
        \renewcommand{\to}{\rightarrow}
        
        \renewcommand{\S}{\Sigma}
        
	\newcommand{\D}{{\cal D}}
        \newcommand{\tr}{{\rm tr}}

        \newcommand{\ad}{{\rm ad}}

        \newcommand{\we}{\wedge}
        \newcommand{\tensor}{\otimes}

        \newcommand{\A}{{\cal A}}
        
        \newcommand{\Comp}{{\rm C}}   
	\newcommand{\R}{{\rm R}}
        
        \newcommand{\GL}{{\rm GL}}

        \newcommand{\qdim}{{\rm qdim}}
        \newcommand{\Vect}{{\bf Vect}}
        \newcommand{\Cob}{{\bf Cob}}
        \newcommand{\C}{{\bf C}}
         \newcommand{\K}{{\bf K}}

        \newcommand{\be}{\begin{equation}}
        \newcommand{\ee}{\end{equation}}
        \newcommand{\ba}{\begin{eqnarray}}
        \newcommand{\ea}{\end{eqnarray}}
      \newcommand{\ban}{\begin{eqnarray*}}
        \newcommand{\ean}{\end{eqnarray*}}
        \newcommand{\barr}{\begin{array}}
      \newcommand{\earr}{\end{array}}

\begin{document}

      \begin{center}
      {\bf 4-Dimensional $BF$ Theory with Cosmological Term\\
      as a Topological Quantum Field Theory\\ }
      \vspace{0.3 cm}
       { John C. Baez\\}
      \vspace{0.3cm}
      {\small Department of Mathematics \\
      University of California\\
        Riverside, CA 92521\\}
     {\small email:  baez@math.ucr.edu}
     \end{center}

\begin{abstract}
Starting from a Lie group $G$ whose Lie algebra is equipped with an
invariant nondegenerate symmetric bilinear form,
we show that 4-dimensional $BF$ theory with cosmological term
gives rise to a TQFT satisfying a generalization of Atiyah's axioms
to manifolds equipped with principal $G$-bundle.   The case $G = \GL(4,\R)$
is especially interesting because every 4-manifold is then naturally
equipped with a principal $G$-bundle, namely its frame bundle.
In this case, the partition function of a compact oriented
4-manifold is the exponential of its signature, and the resulting
TQFT is isomorphic to that constructed by Crane and Yetter
using a state sum model, or by Broda using a surgery presentation
of 4-manifolds.
\end{abstract}

\section{Introduction}
In comparison to the situation in 3 dimensions, topological quantum field
theories (TQFTs) in 4 dimensions are poorly understood.   This
is ironic, because the subject was initiated by an attempt to understand
Donaldson theory in terms of a quantum field theory in 4 dimensions.
However, this theory has never
been shown to fit Atiyah's \cite{Atiyah2} subsequent axiomatic description
of a TQFT, and it is unclear whether one
should even expect it to.

Here we consider a much simpler theory, $BF$ theory with cosmological term
\cite{Horowitz,BT2,CCFM,CM}.
This theory depends on a choice of Lie group $G$ whose Lie algebra
is equipped with an invariant nondegenerate symmetric bilinear
form.   Given an oriented 4-manifold $M$ equipped with principal
$G$-bundle $P$, the fields in this theory are a connection $A$ on $P$
together with an $\ad P$-valued
2-form $B$ on $M$.  The Lagrangian is given by
\[       \tr (B \we F + {\Lambda\over 12} B \we B) ,\]
where it is crucial for our purposes that
$\Lambda$ be nonzero.   In certain cases
this theory can be regarded as a simplified version of general
relativity \cite{Baez}.  This the origin of the term
`cosmological constant' for $\Lambda$ and the curious factor of ${1\over 12}$.
Indeed, this one of our main reasons for studying the theory,
but we defer further discussion of this aspect to another paper.

Starting with this Lagrangian and performing some nonrigorous
computations, we obtain results which we then take as the
definition of a TQFT-like structure satisfying the obvious generalization of
Atiyah's axioms to the case of manifolds equipped with principal $G$-bundle.
Choosing $G = \GL(4,\R)$, we then obtain a TQFT satisfying Atiyah's
axioms by letting $P$ be the frame bundle of $M$.   As it turns out,
the partition function of any compact oriented 4-manifold $M$ is then
$\exp(-36 \pi^2 i\sigma(M)/\Lambda)$,
where $\sigma(M)$ is the signature of $M$.  This fact
says the theory is uninteresting as far as new 4-manifold invariants are
concerned.  However, it is the key to proving a conjecture
that has been discussed in the mathematical physics community for
some time: namely, that the Crane-Yetter-Broda
theory is isomorphic to $BF$ theory.

To understand the origins of the Crane-Yetter-Broda theory and this conjecture
about it, note that one may also define $BF$ theory with cosmological term in
dimension 3.  This has the Lagrangian
\[          \tr(B \we F + {\Lambda \over 3} B \we B \we B) ,\]
where now $B$ is an $\ad P$-valued 1-form.
Although it is difficult to state a concise theorem to this effect,
it is by now commonly accepted that the quantum version of
this theory is isomorphic as a TQFT to a state sum model of
Turaev-Viro type \cite{CCFM,Roberts1,T,TV,Walker,Witten1,Witten2}.
The original Turaev-Viro model was defined using the quantum
group $U_q\Sl(2)$, but it was subsequently generalized
to other quantum groups \cite{BW,Yetter}, and one expects 3-dimensional
$BF$ theory for any simply-connected compact semisimple group $G$ to
be isomorphic to the state sum model based on the corresponding
quantum group $U_q \Comp \g$, with the value of $q$ depending on $\Lambda$.

Using a 4-dimensional state sum model similar to that of Turaev and Viro,
Crane and Yetter \cite{CY}
succeeded in obtaining a 4-dimensional TQFT from $U_q \Sl(2)$.
Shortly thereafter Broda \cite{Broda} constructed a similar theory using a
surgery
presentation of 4-manifolds.   Then Roberts \cite{Roberts1} showed that the
Crane-Yetter and the Broda theories were essentially the same, and that when
properly normalized this TQFT gives as the
partition function of any compact 4-manifold just the exponential of its
signature.  All these results have been extended to other quantum groups by
Crane, Kauffman, and Yetter \cite{CKY2}.

While not yielding new 4-manifold invariants, the Crane-Yetter-Broda
theory is still rather interesting.
First, there is a close relationship between this theory
and Chern-Simons theory.  For example, Roberts \cite{Roberts3}
has shown that at least for $U_q \Sl(2)$,
if one takes a triangulated compact oriented 4-manifold $M$ with boundary
$\Sigma$ and computes the state sum over all labellings, not
holding the labels fixed on $\Sigma$, one obtains the Chern-Simons partition
function
$Z_{CS}(\Sigma)$.  This result probably holds quite generally for other
quantum groups as well.
This is especially nice because no 3-dimensional state sum
model for Chern-Simons theory is known.  In particular, the
phase of $Z_{CS}(\Sigma)$ depends on a choice of framing for $\Sigma$
in a manner that seems difficult to incorporate into
a 3-dimensional state sum, but in the 4-dimensional state sum
the framing is determined by the choice of  bounding
4-manifold $M$.    Second, it is interesting to have a simple
state-sum formula for the signature.  As noted by Crane, Kauffman,
and Yetter \cite{CKY2}, this ``allows us to factor 4-manifold signatures
along any 3-manifold.''  These two points are closely related,
because two 4-manifolds with boundary equal to $\Sigma$ determine the
same framing of $\Sigma$ if and only if their signatures are equal.

For various reasons, it became natural
to suspect that the Crane-Yetter theory corresponds to $BF$ theory in dimension
4 in much the same way that the Turaev-Viro theory corresponds to $BF$ theory
in dimension 3.  By giving a purely differential-geometric construction
of 4-dimensional $BF$ theory with $G = \GL(4,\R)$ as a TQFT,
and proving that this TQFT is isomorphic to the Crane-Yetter theory, we
obtain a precise result to this effect.

We hope this result sheds some
light on the correspondence between two important but quite
different approaches to topological quantum
field theory: the differential-geometric approach and the combinatorial
approach.  However, many details of this correspondence remain to be
understood.  We discuss some particular directions for further
exploration in the Conclusions.

\section{$BF$ Theory}

In this section we first review some basic properties of $BF$ theory
with cosmological term in 4 dimensions,
already established by heuristic arguments due to Blau and Thompson
\cite{BT2}, Horowitz \cite{Horowitz}, Cattaneo, Cotta-Ramusino, Fr\"ohlich
and Martellini \cite{CCFM,CM} and others.
Then we rigorously construct a 4-dimensional TQFT having
these properties, and show it satisfies a generalization of the
Atiyah axioms.

Let $G$ be a Lie group, either
real or complex, with an invariant nondegenerate symmetric
bilinear form on its Lie algebra $\g$.
Let $M$ be an oriented 4-manifold equipped with a principal $G$-bundle $P_M$
over it.   Let $A$ be a connection on $P$ and $B$ an $\ad P$-valued 2-form.
Then the 4-dimensional $BF$ action with
cosmological constant $\Lambda\ne 0$ is given by
\be    S_{BF}(A,B) = \int_M \tr(B \we F + {\Lambda\over 12} B \we B).
\label{SBF4} \ee
Here the `trace' denotes the use of the bilinear form to
turn the quantities in parentheses into an ordinary 4-form.

Heuristically we expect that this action gives rise to
a TQFT $Z_{BF}$ having a one-dimensional space
of states for any compact oriented 3-manifold $\S$.  Arguments
for this have been given using both canonical and
and path-integral approaches.
In the canonical approach, the
kinematical phase space associated to theory on $\R \times \S$ is
the cotangent bundle of the space $\A_\S$ of
connections on $P|_\S$, with the quantity canonically conjugate
to $A|_\S$ being $B|_\S$:
\[            \{B_{ij}^a(x),A_{kb}(y)\} =
 \delta^a_b \epsilon_{ijk} \delta(x,y), \]
where we use spacelike indices $i,j,k,\dots$ and internal indices
$a,b,c,\dots$, and we raise and lower internal indices using the bilinear
form on $\g$.  However, there are constraints: the
Gauss law and the constraint
\ba    F_{ij}^a + {\Lambda \over 6}  B_{ij}^a =0 .\label{con} \ea
Thus in the Dirac approach to quantization,
physical states $\psi$ are functions on $\A_\S$ invariant
under small gauge transformations and satisfying
\be (F_{ij}^a -i {\Lambda \over 6} \epsilon_{ijk}
{\delta \over \delta A_{ka}})\psi =
0 \label{qcon} \ee
Now suppose that $P|_\S$ is trivializable.  Then
this first-order linear PDE has one solution up to a constant
factor, namely
\be         \psi(A) = e^{-{3i\over \Lambda} S_{CS}(A)} ,\label{state} \ee
since the Chern-Simons action
\be   S_{CS}(A) = \int_\Sigma \tr(A \we dA + {2\over 3} A \we A \we A)
\label{SCS} \ee
satisfies
\[  {\delta S_{CS} \over \delta A_{ka}} = \epsilon^{ijk} F_{ik}^a.\]
Note that the definition of the Chern-Simons action, and hence
the state $\psi$, depends on a choice of trivialization of $P|_\S$.
As well-known, the Chern-Simons
action is invariant under small gauge transformations, which we may
interpret as saying that $\psi$ satisfies the Gauss law.   On the other hand,
the Chern-Simons action changes by a constant under large
gauge transformations, which changes $\psi$ by a phase.  This is
not a problem in the present context, since we only expect
(\ref{qcon}) to have a unique solution up to a numerical factor.

The path-integral approach is less rigorous, but it clarifies
the role of the 2nd Chern class, and suggests the ideas
necessary for rigorously constructing a 4d TQFT.
If $M$ has (possibly empty) boundary $\partial M = \S$, we
expect to obtain a vector
$\psi$ in the space of states on $\S$ as follows:
\[ \psi(A_\S) = \int_{A|_\S = A_\S}
\D A \int \D B \;e^{i\int_M \tr(B \we F + {\Lambda\over 12}
B \we B)} .\]
To compute this,
we can first complete the square and do the integral over
$B$:
\[  \psi(A_\S) \propto \int_{A|_\S = A_\S}
\D A  e^{-{3i\over \Lambda}\int_M \tr(F \we F)}.\]
Now suppose again that $P|_\S$ is trivializable.  A choice
of trivialization defines a flat connection on $P|_\S$,
which we extend arbitrarily to a connection $A_0$ on all of
$P$, allowing us to express any other connection on $P$
in terms of an $\ad P$-valued 1-form.  Then
basic relation between the 2nd Chern form and the Chern-Simons form
yields
\be   \int_M \tr(F \we F) = S_{CS}(A_\S) + \int_M \tr(F_0 \we F_0)
\label{reln} \ee
where $F_0$ is the curvature of $A_0$. If we allow ourselves to
neglect the volume factor due to the (heuristic) integrals over
$A$ and $B$, it follows that
\be  \psi(A_\S) = e^{-{3i\over \Lambda}( S_{CS}(A_\S) +
\int_M \tr(F_0 \we F_0)}.\label{reln2} \ee
As expected, $\psi$ is proportional to the solution given
in (\ref{state}).  But also, we find that if the boundary of $M$ is
empty, its partition function equals
\[            Z_{BF}(M) = e^{-{3i\over\Lambda}\int_M \tr(F \we F)}  \]
where $F$ is the curvature of an arbitrary connection
on $P$.  Note that this partition function really depends not
only on $M$ but on the bundle $P$.

While the arguments above are heuristic, we can take the results and use them
to {\it define} a TQFT-like structure satisfying a generalization of Atiyah's
axioms.   Roughly speaking, this will be a functor $Z_{BF}$
from the category $\C$ of `cobordisms between compact
oriented 3-manifolds equipped
with trivializable principal $G$-bundle' to the category $\Vect$ of
vector spaces.
More precisely, an object $\Sigma$ in $\C$ is a compact oriented
3-manifold, also denoted $\Sigma$, together with a trivializable principal
$G$-bundle $P_\Sigma \to \Sigma$.
A morphism $M \maps \S \to \S'$ in $\C$ is an equivalence
class of compact oriented 4-manifolds $M$ with boundary,
equipped with principal $G$-bundle $P_M \to M$ and bundle
isomorphism $\tilde f_M \maps  P_{\overline\S}
 \cup P_{\S'} \to P_M|_{\partial M}$
lifting an orientation-preserving diffeomorphism $f_M \maps
\overline \S \cup \S' \to \partial M $.  The equivalence relation
is that $M \sim M' $ if there is a bundle
isomorphism $\alpha \maps P_M \to P_{M'}$ such that $\tilde f_{M'} =
\alpha \circ \tilde f_M$.  We will not always be so pedantic, however:
usually we will work with representatives rather than equivalence classes.

We define the functor $Z_{BF}$ as follows.
If $\S$ is an object in $\C$, let $\A_\S$ denote the
space of connections on $P_\S$, and let
$Z_{BF}(\Sigma)$ be the space of functions on $\A_\S$
that are multiples of $\exp(-{3i\over \Lambda}S_{CS}(A))$.
If $\Sigma$ is empty we set $Z_{BF}(\Sigma) = \Comp$.
While we need a trivialization of $P_\S$ to define the Chern-Simons
action, and the action may change by a constant as we change
our choice the trivialization, the space $Z_{BF}(\S)$ is independent
of this choice.

The spaces $Z_{BF}(\S)$ have some properties
one expects in a TQFT.  First, $Z_{BF}(\S \cup \S') =
Z_{BF}(\S) \tensor Z_{BF}(\S')$.  More precisely, if we
choose any trivialization of $P_{\S} \cup P_{\S'}$, and restrict
this to trivializations of $P_\S$ and $P_{\S'}$, we have
\[       \exp(-{3i\over \Lambda}S_{CS}(A)) =
   \exp(-{3i\over \Lambda}S_{CS}(A|_{\S}))\,
   \exp(-{3i\over \Lambda}S_{CS}(A|_{\S'}))  \]
for any connection $A$ on $P_\S \cup P_{\S'}$.
This gives an isomorphism $Z_{BF}(\S \cup \S') \simeq
Z_{BF}(\S) \tensor Z_{BF}(\S')$, which one can check does not
depend on the choice of trivialization.  Second, if $\overline \S$ denotes
$\S$ with its orientation reversed, but with the same bundle $P_\S$ over it,
then
$Z_{BF}(\overline \Sigma) = Z_{BF}(\Sigma)^\ast$.  More
precisely, reversing the orientation
switches the sign of the Chern-Simons action, so if
$\psi \in Z_{BF}(\overline \Sigma)$ and $\phi
\in Z_{BF}(\Sigma)$, the product $\psi \phi$ is a constant
function on $\A_\S$, which can be identified with a number times the
constant function $1$.   This gives an isomorphism
$Z_{BF}(\overline \Sigma) \simeq Z_{BF}(\Sigma)^\ast$, which
again does not depend on the choice of trivialization.

If  $M \maps \emptyset \to \S$ is a morphism in $\C$,
$Z_{BF}(M)$ should be a linear map from the complex
numbers to $Z_{BF}(\S)$, or equivalently,
a vector $\psi \in Z_{BF}(\S)$.  We define this vector
by
\be    \psi(A_\S) = e^{-{3i\over \Lambda}
\int_M \tr(F \we F)},\label{zbf} \ee
where $F$ is the curvature of any connection $A$ extending
$A_\S$ to all of $P_M$.  By equation (\ref{reln}),
$\psi$ lies in $Z_{BF}(\S)$ and is independent of the choice of
$A$.

More generally, given any morphism $M \maps \S \to \S'$ in $\C$,
since the boundary of $M$ is $\overline \S \cup \S'$, we
obtain a vector in $Z_{BF}(\S)^\ast \tensor Z_{BF}(\S')$ by
the above procedure.  We then define $Z_{BF}(M) \maps
Z_{BF}(\S) \to Z_{BF}(\S')$ to be this vector, reinterpreted
as a linear map from $Z_{BF}(\S)$ to $Z_{BF}(\S')$.

It is easy to see that $\C$ becomes a rigid symmetric monoidal category
in a manner analogous to the usual categories of cobordisms, with
the tensor product of objects and morphisms being
given by disjoint union, and duality for objects
being given by orientation reversal.  We then have:

\begin{theorem}\et  $Z_{BF} \maps \C \to \Vect$ is a
symmetric monoidal functor.     \end{theorem}

Proof -   A straightforward computation.
\qed

\noindent Using a simple argument noticed by Crane and
Yetter \cite{CY2}, it follows that $Z_{BF}$ preserves duals up to
canonical isomorphism.

\section{The Case $G = \GL(4,\R)$}

The construction of the previous section could in fact be generalized
to higher even dimensions using other characteristic classes and their
secondary characteristic classes.  Just as in
Theorem 1, these would give functors from a category
$\C$ of `cobordisms between
compact oriented $(n-1)$-manifolds equipped with trivializable
principal $G$-bundle' to the category $\Vect$.
What makes the construction particularly interesting when $n = 4$ is that
the tangent bundle of any compact oriented 3-manifold is trivializable.
This yields various procedures for obtaining objects and
morphisms in $\C$ from those in $4\Cob$, the category
of cobordisms between compact oriented 3-manifolds.
A procedure of this sort that involved
no arbitrary choices might yield a functor
from $4\Cob$ to $\C$, which we could compose with $Z_{BF}$ to obtain
a functor from $4\Cob$ to $\Vect$, and thus, with some luck,
a TQFT.   In fact, the procedure we consider gives a map from
$4\Cob$ to $\C$ which is
not quite a functor, but for which the composite with $Z_{BF}$ is
still a TQFT.   This is the $BF$ theory naturally associated to the
frame bundle of a 4-manifold.

Let $G = \GL(4,\R)$, and equip its Lie algebra with the
bilinear form
\[     \langle S,T \rangle = \tr(ST)  ,\]
the trace being that of $4 \times 4$ matrices.
Given any compact oriented 4-manifold with boundary
$M$, the oriented frame bundle $P_M$ is a principal $\GL(4,\R)$-bundle
over $M$.   We can also construct an $\GL(4,\R)$-bundle
over a compact oriented 3-manifold $\S$ as follows.  Let
$T\S$ be its tangent bundle and let $L\S$ be the
trivial line bundle $\S \times \R$ over $\S$.  Then the bundle $P_\S$ of
oriented frames of $T\S \oplus L\S$ is a principal $\GL(4,\R)$-bundle.
Since $T\S$ is trivializable, so is $P_\S$, so we have a way to get
objects in $\C$ from objects in $4\Cob$.  We do not, however, have a systematic
way to get morphisms in $\C$ from morphisms in $4\Cob$.
We may think of a morphism $M\maps \S \to \S'$ in $4\Cob$ as a compact
oriented 4-manifold $M$ with boundary, equipped with an orientation-preserving
diffeomorphism $f_M \maps \overline \S \cup \S' \to \partial M$.
(Actually, just as in $\C$, morphisms in
$4\Cob$ are really certain equivalence classes \cite{Sawin}, but
we shall work with representatives and leave the reader to check
that our constructions make sense at the level of equivalence classes.)
To obtain a morphism in $\C$ from this morphism in $4\Cob$, we need to pick a
bundle isomorphism $\tilde f_M \maps P_{\overline \S} \cup P_{\S'} \to
P_M|_{\partial M}$ lifting $f_M$.   There appears to be no way
to do this without an arbitrary choice.  Without loss of generality
we can assume $\overline\S \cup \S' = \partial M$ and that $f_M$ is the
inclusion map.
Then we can obtain $\tilde f_M$ from an orientation-preserving isomorphism
\[                      \eta \maps T\S \oplus L\S \to TM|_\S \]
together with a similar isomorphism for $\S'$.   On $T\S$ we
define $\eta$ to be the inclusion
\[                     T\S \hookrightarrow TM|_\S  ,\]
and defining $\eta$ on $L\S$ amounts to choosing a
section $v$ of $TM|_\S$.  For $\eta$ to be an isomorphism, it is
necessary and sufficient that $v$ be nowhere tangent to
$\S$.  For it to be orientation-preserving, $v$ must be inwards-pointing.
The same holds for $\S'$ except that $v$ must be outwards-pointing.

So in short, this procedure does not quite give a
functor from $4\Cob$ to $\C$, but only a functor from a
category $4\Cob'$ in which the morphisms $M \maps \S \to \S'$ in
$4\Cob$ are equipped with a bit of extra structure:
a section $v$ of $TM$ over $\partial M$ that is nowhere tangent
to $\partial M$, inwards-pointing on $f_M \S$, and outwards-pointing
on $f_M \S'$.    This should not be surprising:
a similar structure, called the `lapse and shift', plays an important role in
general relativity.  Let us denote this functor from
$4\Cob'$ to $\C$ as $F$.

We then obtain an actual TQFT as follows.  First we define a map
$G \maps 4\Cob \to 4\Cob'$, taking objects to objects and morphisms
to morphisms, but not a functor, as follows.  For each object $\S$ of
$4\Cob$, let $G(\S) = \S$.  For each morphism $M\maps \S \to \S'$ of
$4\Cob$, let $G(M)$ be $M$ equipped with an arbitrary section $v$ of
$TM$ over $\partial M$ with the necessary properties.
Then let $Z \maps 4\Cob \to \Vect$ be the composite $Z_{BF} \circ F
\circ G$.

\begin{theorem}\et $Z$ is a TQFT, that is, a symmetric
monoidal functor from $4\Cob$ to $\Vect$.  \end{theorem}

First we show that $Z$ is a functor.  Note that the only way $G$ fails
to be a functor is that it fails to send identity morphisms to identity
morphisms, so we only need to check that $Z$ has this property.  Let
$\S$ be an object of $4\Cob$ and $M \maps \S \to \S$ the identity on
$\S$.   $Z(M) \maps Z(\S) \to Z(\S)$ may
be identified with a vector $\psi \in Z(\overline \S \cup \S)$.
Let us fix, once and for all, a trivialization of $T\S$.   From this
we obtain a trivialization of $P_\S$, and, using the
standard diffeomorphism $\S \simeq \overline \S$,
also a trivialization of $P_{\overline\S}$.  Then to show that
$Z(M)$ is the identity it suffices to check that
\be \psi(A_{\overline \S \cup \S}) =
e^{-{3i\over \Lambda} S_{CS}(A_{\overline \S \cup \S})}  \label{check}  \ee
for any connection $A_{\overline \S \cup \S}$ on $P_{\overline \S} \cup P_\S$.

For the purposes of checking this, let us identify the
manifold $M$ with $[0,1] \times \S$.
To construct $F(M)$ we have arbitrarily equipped $M$
with a section $v$ of $TM$ over $\partial M$ that is nowhere tangent to
$\partial M = \{0,1\} \times \S$, inwards-pointing on $\{0\} \times \S$,
and outwards-pointing on $\{1\} \times \S$.
This gives an isomorphism $\tilde f_M \maps P_{\overline\S} \cup P_{\S} \to
P_M|_{\partial M}$.   By equation (\ref{zbf}), $\psi$ is given by
\[     \psi(\tilde f_M^\ast A) = e^{-{3i\over \Lambda} \int_M \tr(F \we F)},\]
where $A$ is any connection on $P_M|_{\partial M}$, and $F$ is its curvature.

We can evaluate the right-hand side of the above equation using
equation (\ref{reln}).   Using our trivialization of $T\S$ and the
standard vector field $\partial_t$ on $M = [0,1] \times S$
associated with the coordinate
$t$ on $[0,1]$, we obtain a trivialization of $TM$, hence of $P_M$.
In equation (\ref{reln}) take $A_0$ to be the flat connection on $P_M$
associated with this trivialization, and restrict this
trivialization to $\partial M$ to define the Chern-Simons action.
Then we obtain:
\[    \psi(\tilde f_M^\ast A) = e^{-{3i\over\Lambda}
S_{CS}(A)} . \]
Equation (\ref{check}), which we wish to check, will follow upon
setting $A_{\overline \S \cup \S} = \tilde f_M^\ast A$ if we can show
\be    S_{CS}(A) = S_{CS}(\tilde f_M^\ast A)  .\label{check2} \ee
The key point is to show that the trivialization of $P_M|_{\partial M}$
used to the define the left-hand side, and the trivialization of $P_{\overline
\S \cup \S}$ used to define the right-hand side, are compatible.  It is
{\it not} true that the isomorphism $\tilde f_M$ carries the trivialization of
$P_{\overline \S \cup \S}$ to the trivialization of $P_M|_{\partial M}$.
However, note $v$ is inwards-pointing on $\{0\} \times S$ and
outwards-pointing on $\{1\} \times s$, and so is $\partial_t$.  Thus
$\tilde f_M$ does carry the trivialization of
$P_{\overline \S \cup \S}$ to that of $P_M|_{\partial M}$ up to
a small gauge transformation.   Since the Chern-Simons action
is invariant under small gauge transformations, equation (\ref{check2})
holds.

It is easy to check that $Z$ is monoidal, and an argument like the above
one shows that $Z$ is symmetric.   \qed

\noindent Again, while the fact that $Z$ preserves duals up to canonical
isomorphism
is often taken as part of the definition of a TQFT, this actually follows
from $Z$ being a symmetric monoidal functor \cite{CY2}.

If $M$ is a compact oriented 4-manifold and $A$ is any
connection on $TM$, the Hirzebruch signature
theorem implies that $\int_M \tr(F \we F) = 12 \pi^2 \sigma(M)$,
where $\sigma(M)$ is the signature of $M$.     Thus we have
\[               Z(M) = e^{-36 \pi^2 i \sigma(M)/Lambda}  .\]

\section{Relation to the Crane-Yetter-Broda Theory}

To establish an isomorphism between the Crane-Yetter-Broda
theory and the TQFT of the previous section, we
do not need to know much about the Crane-Yetter-Broda theory.
All we need, in fact,
is that it is a 4-dimensional TQFT such that the partition function of
any 4-manifold equals $\exp(\alpha \sigma(M))$ for some constant
$\alpha$; as we shall see, there is
a unique such TQFT for any $\alpha$.
Crane and Yetter's description of the
theory in terms of a state sum model is interesting nonetheless,
because the existence of such a model means that we have an {\it extended}
TQFT in the sense of Lawrence \cite{Lawrence}.
Thus we begin with a brief review of this state sum model.

The original Crane-Yetter theory was constructed using the representation
theory of the quantum group $U_q\Sl(2)$ when $q$ is a suitable
root of unity.  This was subsequently generalized by Crane, Kauffman and
Yetter \cite{CKY2} to other quantum groups, and is this generalization
that we describe here.  The input to the theory is an semisimple tortile
tensor category $\K$ with trivial center.   The reader may turn to the above
reference for the definitions involved, but the main examples to keep in
mind are certain subquotients of the categories of finite-dimensional
representations of the quantum groups $U_q(\g)$, where
$\g$  is a complex semsimple Lie algebra and $q$ is a suitable root of
unity \cite{Sawin}.
These subquotients are formed by first considering the full subcategory of
completely reducible representations, and then quotienting by the tensor
ideal generated by representations of quantum dimension $0$.

Being semisimple, the category $\K$ has a basis $S$ of simple objects.
We choose this basis
so that it contains the unit $1$ for the tensor product in $\K$, and so
that $a \in S$ implies $a^\ast \in S$.   For any objects $a,b,c \in S$, we
choose a basis $B^{ab}_c$ of the vector space ${\rm Hom}_{\K}(a \tensor b,c)$,
and we write $B$ for the disjoint union of all these bases.
In the quantum group case the objects in $S$ are in 1-1 correspondence with
irreducible representations of nonzero quantum dimension, and we may
think of $B^{ab}_c$ as a basis of intertwining operators from $a \tensor b$
to $c$, modulo intertwining operators that factor through a representation
of quantum dimension zero.  In general, every object $x$ of $\K$ has a
`quantum dimension' $\qdim(x)$, which is typically not an integer.

Let $M$ be a triangulated compact oriented
4-manifold with an ordering of its vertices.  Let
$T_i$ denote the set of (nondegenerate) $i$-simplices in $M$, and let
$n_i = |T_i|$.  By
a `coloring' of $M$, we mean maps
\[  \lambda \maps T_2 \cup T_3 \to S, \qquad \lambda^+,
\lambda^- \maps T_3 \to B \]
such that
\[    \lambda^+(a,b,c,d) \in
B^{\lambda(b,c,d) \, \lambda(a,b,d)}_{\lambda(a,b,c,d)},
\]
and
\[   \lambda^-(a,b,c,d) \in
B^{\lambda(a,c,d) \, \lambda(a,b,c)}_{\lambda(a,b,c,d)},\]
where $a,b,c,d$ are vertices with $a < b < c < d$, and triples and quadruples
of vertices denote triangles and tetrahedra, respectively.
The reader can visualize such a coloring as follows.
We form a graph in the 3-skeleton of $M$ whose intersection with any
tetrahedron
$\tau = (a,b,c,d) \in T_3$ appears as in Figure 1.  The 4 edges of this
graph intersecting
the faces of $\tau$ are labelled with objects in $S$ using $\lambda|_{T_2}$,
while the edge in the center of $\tau$ is labelled with an object
using $\lambda|_{T_3}$.  The 2 vertices are labelled with elements of
$B$ using $\lambda^+$ and $\lambda^-$.

The Crane-Yetter invariant $Z_{CY}(M)$ is then given by
\be   \sum_{{\rm colorings}}
N^{n_0 - n_1} \prod_{\sigma \in T_2} \qdim(\lambda(\sigma))
\prod_{\tau \in T_3} \qdim(\lambda(\tau))^{-1}  \prod_{\rho \in T_4}
|\lambda,\rho|   ,\label{statesum} \ee
where
\[ N = \sum_{a \in S} \qdim(a)^2 , \]
and $|\lambda, \rho|$ is a quantity obtained using the standard
Reshetikhin-Turaev procedure \cite{RT}
from the portion of the graph contained in
the boundary of the 4-simplex $\rho$.

The quantity $Z_{CY}(M)$ turns out to be independent of the triangulation of
$M$.   As shown by Roberts \cite{Roberts1} in the $U_q\Sl(2)$ case and by
Crane, Kauffman and Yetter in general,
\[    Z_{CY}(M) = N^{\chi(M)/2} y^{\sigma(M)}  \]
where $y = (a_+/a_-)^{1/2}$, with $a_\pm$ obtained using
the Reshetikhin-Turaev procedure from the $\pm 1$-framed unknot labelled
with the linear combination of objects
\[   \Omega = \sum_{a \in S} \qdim(a) a .\]
Of course, by multiplying the state sum (\ref{statesum})
by an appropriate exponential of Euler characteristic
$\chi(M) = n_0 - n_1 + n_2 - n_3 + n_4$,
we can easily obtain a modified Crane-Yetter state sum for which the
partition
function of $M$ contains any desired exponential of $\chi(M)$.  In fact,
Broda's \cite{Broda}
approach using a surgery presentation of $M$ very naturally corresponds
to the
modified state sum in which the partition function is just
\[       Z'(M) = y^{\sigma(M)}  .\]
In what follows we use this normalization.

So far we have discussed the theory only for compact 4-manifolds.
One may extend the theory to a TQFT using the following
procedure.   Let $4\Cob_{\rm PL}$
be the piecewise-linear analog of the category $4\Cob$.  For each
object $\S$ of $4\Cob_{\rm PL}$, define $Z'(\S)$ to be the vector space
having as a basis all the morphisms $M \maps \emptyset \to \S$, modulo
those linear combinations $\sum c_M M$ such that
\[              \sum c_M Z'(NM) = 0    \]
for all morphisms $N \maps \S \to \emptyset$.  Then, for each morphism
$N \maps \S \to \S'$, define $Z'(N)$ by
\[              Z'(N) [\sum c_M M] = [\sum c_M NM] .\]
Adapting the ideas of
Blanch\'et {\it et al} \cite{BHMV} to the case of TQFTs that
are not necessarily unitary, one can check that
$Z'$ is well-defined on morphisms and actually gives a TQFT.

Now, while the Crane-Yetter state sum is defined in the piecewise-linear
category, the categories $4\Cob$ and $4\Cob_{\rm PL}$ are equivalent,
so we can transfer the resulting TQFT to one in the smooth category.  It
will cause no confusion in what follows to also call this smooth
version of the Crane-Yetter-Broda theory $Z'$.
Our main result is then:

\begin{theorem}\et  If  $y = \exp(-36 \pi^2 i/Lambda)$,
then $Z \maps 4\Cob \to \Vect$ of Theorem 2 and the
Crane-Yetter-Broda theory $Z' \maps 4\Cob \to \Vect$ are equivalent
as TQFTs.  More precisely, there is a monoidal natural
isomorphism $F \maps Z \to Z'$.  \end{theorem}

Proof - In fact, we shall show that any two TQFTs $Z, Z' \maps
4\Cob \to \Vect$ having
\[           Z(M) = Z'(M) = y^{\sigma(M)} \]
for all compact oriented 4-manifolds $M$ are equivalent
in this sense, where $y \ne 0$.  Concretely, this
means, first of all,
that for any object $\S$ of $4\Cob$ there is an isomorphism
$F_\S \maps Z(\S) \to Z'(\S)$.   This isomorphism should be
natural, in the sense that for any morphism
$M \maps \S_1 \to \S_2$, one has $Z'(M)F_{\S_1} = F_{\S_2}Z(M)$.
It should also be compatible with the monoidal structure.
This compatibility condition is easiest to state if
we use Mac Lane's theorem to replace $4\Cob$ and $\Vect$ by
equivalent strict monoidal categories; it then states that
$F_{\S_1 \cup \S_2} = F_{\S_1} \tensor F_{\S_2}$
and $F_{\emptyset} = {\rm id}$.

To construct an isomorphism with these properties,
note first that for any object $\S$ of $4\Cob$, $Z(\S)$ is
1-dimensional, and similarly for $Z'(\S)$.  To see this,
recall that
\[          \dim Z(\S) = Z(S^1 \times \S) ,\]
and that $\sigma(S^1 \times \S) = 0$.
Also note that for any object $\S$ of $4\Cob$ there is
a morphism $M \maps \emptyset \to \S$, because the oriented
cobordism group vanishes in dimension 3.  Moreover,
for any morphism $M \maps \emptyset \to \S$ the vector
$Z(M)1$ is nonzero, because $\overline M$ is a morphism from
$\S$ to $\emptyset$, and
\[    Z(\overline M)Z(M)1 = Z(\overline M M)1 = y^{\sigma(\overline MM)}
\ne 0.  \]

We claim that for any $\S$, there is a unique isomorphism $F_\S \maps
Z(\S) \to Z'(\S)$ such that $F_\S(Z(M)1) = Z'(M)1$ for all $M \maps
\emptyset \to \S$.  Uniqueness is clear, because the remarks above
imply $Z(\S)$ is spanned by $Z(M)1$.  For existence, we need to check
that if $M_1, M_2 \maps \emptyset \to \S$, so that $Z(M_1) = \alpha Z(M_2)$
and $Z'(M_1) = \alpha' Z'(M_2)$, then $\alpha = \alpha'$.
To see this it suffices to note that
\[         Z(\overline M_1 M_1) = \alpha Z(\overline M_1 M_2),\qquad
Z'(\overline M_1 M_1) = \alpha' Z'(\overline M_1 M_2), \]
and that $Z$ and $Z'$ agree and are nonzero on compact 4-manifolds.

To show that $F_\S$ is natural it suffices to consider a morphism
$M \maps \S_1 \to \S_2$ and check that $Z'(M)F_{\S_1} = F_{\S_2}Z(M)$
on a vector of the form $Z(N)1$, where $Z(N)\maps \emptyset \to \S_1$.
This is clear:
\ban        Z'(M)F_{\S_1}Z(N)1 &=&  Z'(M)Z'(N)1 \\
&=& Z'(MN) 1 \\
&=& F_{\S_2} Z(MN) 1 \\
&=& F_{\S_2} Z(M) Z(N) 1. \ean

To show that $F_\S$ is compatible with the monoidal structure
choose $M_1 \maps \emptyset \to \S_1$ and $M_2 \maps \emptyset \to \S_2$.
Then we have
\ban       F_{\S_1 \cup \S_2} Z(M_1 \cup M_2) 1 &=& Z'(\S_1 \cup \S_2)1 \\
&=&   Z'(\S_1)1 \tensor Z'(\S_2)1 \\
&=&    F_{\S_1}Z(\S_1)1 \tensor F_{\S_2}Z(\S_2)1 \\
&=&    (F_{\S_1} \tensor F_{\S_2}) Z(M_1 \cup M_2)1 \ean
so $F_{\S_1 \cup \S_2} = F_{\S_1} \tensor F_{\S_2}$.  Similar manipulations
show that $F_\emptyset$ is the identity.
\qed

\section{Conclusions}

We have exhibited 4-dimensional $BF$ theory with cosmological
term and quite general gauge group $G$
as a TQFT-like functor, and obtained an actual TQFT
by setting $G = \GL(4,\R)$ and using the frame
bundle as a natural choice of $G$-bundle.   This
TQFT is equivalent to the Crane-Yetter-Broda theory when the
the cosmological constant $\Lambda$ and the constant $y$ appearing in
the Crane-Yetter-Broda theory are related as in Theorem 3.
However, the real reason for this equivalence
is still somewhat mysterious, and deserves more study.

Indeed, it might seem surprising at first that
the Crane-Yetter-Broda theory depends so little on the choice of
category $\C$ --- or, more concretely, on the choice of
quantum group $U_q \g$.  The labels involved in the
Crane-Yetter state sum depend in detail upon this choice, so
different choices should give different extended TQFTs ---
that is, they should assign different algebraic data to manifolds
with corners, such as the 4-simplex itself.
However, regarded just as a TQFT the Crane-Yetter-Broda theory
depends on $U_q \g$ only through the single constant $y$!

To understand this, note that
in quantum group examples $y$ is the quantity by which the Chern-Simons
partition function
of a compact oriented framed 3-manifold is multiplied
when one changes the framing by one `twist' \cite{RT}.
Thus one may say that as a TQFT
the Crane-Yetter-Broda theory simply keeps track of the
framing-dependence of the corresponding Chern-Simons theory.
As an extended TQFT, however, it appears to contain enough information
to compute the Chern-Simons partition function --- though so far this
has been proved only for $U_q \Sl(2)$.
Different extended TQFTs may correspond to the same TQFT, as is
already clear from 2-dimensional examples, so the fact that we have
found a $BF$ theory isomorphic to the Crane-Yetter-Broda theory
as a TQFT does {\it not} imply that this $BF$ theory
is isomorphic to it as an extended TQFT.  This is relevant to
the problem of obtaining invariants of embedded surfaces
from $BF$ theory \cite{CCFM,CM}, since obtaining such invariants
uses its structure as an extended TQFT.  It is known how to
obtain such invariants from the Crane-Yetter-Broda
theory \cite{Roberts2,Yetter2}, but until we obtain an isomorphism
between this theory and a $BF$ theory as extended TQFTs,
we will not know what these results imply about $BF$ theory.

Now, just as the Crane-Yetter-Broda theory as an extended
TQFT appears to be able to compute the corresponding Chern-Simons
partition function, equation (\ref{reln2}) suggests a similar
result for $BF$ theories.  So there is reason to conjecture that a suitable
$BF$ theory is isomorphic to the Crane-Yetter-Broda theory
as an extended TQFT.  In fact, it seems that for $G$ simply-connected,
compact and semisimple, the Crane-Yetter-Broda
theory with quantum group $U_q \Comp \g$ corresponds to a $BF$ theory
with gauge group $G \times \GL(4,\R)$, with each 4-manifold being
equipped with a principal bundle given by the trivial $G$-bundle times its
frame bundle.  This is already implicit in Witten's original paper
\cite{Witten2} on Chern-Simons theory, since in his approach the
framing-dependence arises by adding to the Chern-Simons Lagrangian with
gauge group $G$ a `gravitational' Chern-Simons term.  We
hope to return to this conjecture in future work.

\subsection*{Acknowledgements}
The author would like to thank David Johnson,
Allen Knutson, Greg Kuperberg, Justin Roberts, Stephen Sawin and Oleg
Viro for helpful discussion and correspondence, and
especially Alberto Cattaneo, Paulo Cotta-Ramusino
and Maurizio \break Martellini for conversations and hospitality
at the University of Milan.

\end{document}